\documentclass[3p,sort&compress,11pt,twoside]{elsarticle}
\usepackage[T1]{fontenc}
\usepackage[utf8]{inputenc}
\usepackage[english]{babel}
\usepackage{amsmath,amssymb}
\usepackage{times}
\usepackage{pgfplots}
\pgfplotsset{compat=newest}
\usepackage{hyperref}
\newcommand{\bbR}{\mathbb{R}}
\newcommand{\Rd}{{\bbR^d}}
\newcommand{\assign}{\leftarrow}
\newcommand{\diff}{\mathrm{d}}

\newcommand{\dx}{\diff x}
\newcommand{\dd}[2][{}]{\frac{\diff#1}{\diff#2}}
\newcommand{\ddt}{\dd[]t}
\newcommand{\bz}{\bar z}
\newcommand{\bZ}{\bar Z}
\newcommand{\hZ}{\widehat Z}
\newcommand{\bzeta}{\bar\zeta}
\newcommand{\hzeta}{\widehat\zeta}
\newcommand{\hj}{{\hat\jmath}}
\DeclareMathOperator{\opDyn}{\mathbf{D}}
\DeclareMathOperator{\opFlip}{\mathbf{F}}
\DeclareMathOperator{\opInt}{\mathbf{I}}
\DeclareMathOperator{\oprand}{\mathbf{r}}
\DeclareMathOperator{\opRand}{\mathbf{R}}
\newcommand{\distrN}{\mathcal{N}}
\newcommand{\codeline}[1]{\texttt{\textbf{#1}}}

\hypersetup{
  unicode=true,
  breaklinks=true,
  colorlinks=true,
  linkcolor=blue,
  citecolor=blue,
  filecolor=blue,
  urlcolor=MidnightBlue,
  pdftitle={Extra Chance Generalized Hybrid Monte Carlo},
  pdfauthor={Cédric M. Campos, J. M. Sanz-Serna},
  pdfsubject={PACS 02.50.Ng, 02.70.Ns, 02.70.Tt, 05.10.Ln, 31.15.xv; MSC2010 65C05, 65C35, 82C80, 74A25},
  pdfkeywords={sampling methods, hybrid Monte Carlo, detailed balance, delayed rejection, rejection avoidance, molecular dynamics},
  pdfcreator={GNU Emacs 24.3.1},
}

\begin{document}
\title{Extra Chance Generalized Hybrid Monte Carlo\tnoteref{xhmc}}
\tnotetext[xhmc]{A MATLAB implementation of the algorithm and other utilitarian scripts can be found at \url{http://github.com/vitaminace33/xhmc}.}
\author{Cédric~M. Campos\corref{cor}}
\ead{cedricmc@\{uva,icmat\}.es}
\address{Dept. Matemática Aplicada e IMUVA, Fac. Ciencias, UVA\\Paseo de Belén 7, 47011 Valladolid, Spain}
\cortext[cor]{Corresponding author}
\author{J.~M. Sanz-Serna\corref{nocor}} 
\ead{jmsanzserna@gmail.com}
\address{Dept. Matemáticas, Universidad Carlos III de Madrid\\Edificio Sabatini, 28911 Leganés (Madrid), Spain}

\begin{abstract}
We study a method, Extra Chance Generalized Hybrid Monte Carlo, to avoid rejections in the Hybrid Monte Carlo method and related algorithms. In the spirit of delayed rejection, whenever a rejection would occur, extra work is done to find a fresh proposal that, hopefully, may be accepted. We present experiments that clearly indicate that the additional work per sample carried out in the extra chance approach clearly pays in terms of the quality of the samples generated.
\end{abstract}


\maketitle

\section{Introduction}
In this paper we study a technique, Extra Chance Generalized Hybrid Monte Carlo (XCGHMC), to avoid rejections in the Hybrid Monte Carlo (HMC) method \cite{DuKePe87} and its variants (e.g. \cite{Ho91,KePe01,IzHa04,AkRe08,FaSaSk14}). In the spirit of what in the statistics literature is called \emph{delayed rejection} \cite{TiMi99,Mi01,GrMi01}, whenever a rejection would occur, additional work is done to find a fresh proposal that, hopefully, may be accepted. Rejections, while essential to ensure that the algorithms sample from the right target probability distribution, contribute to an increase of the correlation of the samples \cite{Pe73,Ti98}. Furthermore, for algorithms with partial momentum refreshments \cite{Ho91,KePe01}, each rejection necessarily requires a flip of the momentum and interferes with the underlying Hamiltonian dynamics. We present experiments that clearly indicate that the additional work per sample carried out in the extra chance approach in order to avoid rejections clearly pays in terms of the quality of the samples generated.

The algorithm studied here is mathematically equivalent to that suggested by Sohl-Dickstein, Mudigonda and DeWeese in \cite{SoMuWe14}. However the actual formulas for the acceptance probability in the present work are different from those in \cite{SoMuWe14}. We believe that the formulas used here provide insight into the probabilities involved. Furthermore we prove that the extra chance algorithm actually satisfies detailed balance/stochastic reversibility; this is surprising, as \cite{SoMuWe14} suggests that the avoidance of rejections comes at the price of the violation of detailed balance. Detailed balance is of course a valuable property in the performance and analysis of Markov Chain Monte Carlo algorithms \cite{RoCa04}. For instance the estimation of the effective sample size or autocorrelation time used in our experiments (taken from \cite{Ge92}) relies on the chain being reversible with respect to the target distribution.

In turn the basic idea in \cite{SoMuWe14} is not essentially different from delayed rejection \cite{TiMi99,Mi01,GrMi01}. However some peculiarities of HMC (momentum flip, modified detailed balance vs.\ standard detailed balance, proposals being obtained via a deterministic flow) make it difficult, or even impossible, to apply the material of \cite{TiMi99,Mi01,GrMi01} to devise or analyze algorithms similar to those considered in \cite{SoMuWe14} or here.

The literature is not lacking in suggestions to avoid rejection/momentum flips in HMC, see e.g.\ \cite{AkRe08,AkBoRe09a,AkBoRe09b,So12,WaPa12}. Comparisons between those techniques and the extra chance approach are not within our scope here. Neither shall we be concerned with comparing HMC with alternative sampling algorithms as done in \cite{CaLeSt07}.

In Section \ref{sec:algorith} we present the extra chance algorithm. Section \ref{sec:acceptance} is devoted to an analysis of the acceptance probabilities. Proof of concept numerical experiments are reported in Section \ref{sec:numerics} and Section \ref{sec:conclusion} concludes. Some more mathematical results are given in Appendices \ref{sec:stationarity}--\ref{sec:equivalence}.

\section{Algorithm} \label{sec:algorith}
The aim is to obtain samples \(x_n\) from a target probability distribution in the \emph{state space} \(\Rd\) with density

\begin{equation}\label{eq:target}
\frac1Z \exp(-\beta V(x)),\qquad Z = \int_\Rd \exp(-\beta V(x))\,\dx.
\end{equation}
The algorithms considered here do not require that \(Z\) be known; they merely need to evaluate \(V\) and its gradient \(\nabla V\). They generate a Markov chain \(x_0\to x_1\to \dots \to x_N\) \cite{RoCa04} that has \eqref{eq:target} as an invariant distribution, in such a way that under suitable ergodic behavior, it is possible to estimate averages \(\langle A\rangle\) with respect to the target by taking means of the values of \(A\) along a realization of the chain:
\begin{equation}\label{eq:average}
\langle A\rangle = \frac1Z\int_\Rd A(x) \exp(-\beta V(x))\,\dx\approx \frac1{N+1}\sum_{n=0}^N A(x_n).
\end{equation}
The format in \eqref{eq:target} implies that the density is everywhere positive. As shown in e.g.\ \cite{FaSaSk14}, it is not difficult to extend the algorithms and analysis to cases where the density vanishes in a subset of the state space.

Regardless of the details of the application in mind, HMC and its variants use the Hamiltonian formalism of classical mechanics. The components of \(x\) are interpreted as generalized co-ordinates describing the configuration of a mechanical system and an auxiliary \(d\)-dimensional vector \(y\) is introduced whose components represent the associated conjugated momenta. We shall use the notations \(z = (x,y)\), \(z^\prime = (x^\prime, y^\prime)\), etc.\ to refer to points in the \emph{phase space} \(\Rd\times\Rd\). If \(M\) is a user-specified, symmetric positive-definite \(d\times d\) mass matrix, the algorithms use the Hamiltonian function (total mechanical energy)
\begin{equation}\label{eq:hamiltonian}
H(z) = \frac12\, y^TM^{-1}y + V(x),
\end{equation}
and the probability distribution in the phase space defined by the unnormalized density
\begin{equation}\label{eq:distribution}
\rho(z) = \exp(-\beta H(z))=\exp(-(\beta/2)\, y^TM^{-1}y)\times\exp(-\beta V(x)).
\end{equation}
The product structure of \(\rho\) implies that \(x\) and \(y\) are stochastically independent; \(x\) is distributed according to the target \eqref{eq:target} and \(y\sim \distrN(0,M)\), i.e. \(y\) is Gaussian with zero mean and covariance matrix \(M\).

The dynamics associated with \eqref{eq:hamiltonian} is given by
\begin{equation}\label{eq:odes}
\ddt\, x = M^{-1}y,\qquad \ddt\, y = -\nabla V(x),
\end{equation}
a system of differential equations whose solution flow exactly preserves the distribution \eqref{eq:distribution}, see e.g.\ \cite{Sa14}. In practice this flow cannot be computed in closed form and one has to resort to numerical approximations; the Störmer-Verlet/leapfrog integrator \cite{SaCa94,HaLuWa06,LeRe04,Sc10} is the method of choice. Figure~\ref{fig:verlet} shows pseudocode for computing the result \(z^\prime = (x^\prime, y^\prime) = \opInt(z)\) of \(L\) integration time-steps of length \(\Delta t >0\) starting from the initial point \(z\). The transformation \(\opInt\), which maps the phase space into itself, is both \emph{volume preserving} (i.e. has unit Jacobian determinant) and \emph{reversible}. Reversibility means that, for each \(z^\prime\), \(\opInt^{-1}(z^\prime)=\opFlip(\opInt(\opFlip(z^\prime)))\), where \(\opFlip\) denotes the \emph{momentum flip transformation:} \(\opFlip(x,y) = (x,-y)\). Thus to recover the initial point \(z = \opInt^{-1}(z^\prime)\) of a numerical integration it is sufficient to flip the momentum, take \(L\) (forward) time-steps and flip again the momentum; this exactly matches a fundamental property of the true solution flow of \eqref{eq:odes}. Note that \(\opFlip\) is volume preserving and, furthermore, \(H(z) = H(\opFlip(z))\).

The XCGHMC algorithm is summarized in Fig.~\ref{fig:main}. The proof that \(\rho\) in \eqref{eq:distribution} is an invariant density of the chain \(z_0\to z_1\to \dots \to z_N\) is given in Appendix \ref{sec:stationarity}. Appendix \ref{sec:balance} shows that, if \(y_0\sim \distrN(0,M)\), the marginal chain \(x_0\to x_1\to \dots \to x_N\) obtained by taking the \(x\) components satisfies detailed balance with respect to the target distribution \eqref{eq:target}. As pointed out in the introduction, the algorithm in Fig.~\ref{fig:main} is mathematically equivalent to that suggested in \cite{SoMuWe14}; see Appendix \ref{sec:equivalence}. It is in order to point out that the material in the appendices relies on \(\opInt\) being volume-preserving and reversible and is otherwise independent of the properties of the velocity Verlet integrator; it is therefore possible to use alternative integrators in Fig.~\ref{fig:main}, provided that they are both volume-preserving and reversible. For integrators tailored to this task the reader is referred to \cite{BlCaSa14,PrLiEa12} and their references.

\begin{figure}
\begin{center}
\ttfamily\bfseries
\begin{tabular}{cl}
00 & given \(z=(x,y)\in \Rd\times\Rd\) \\
01 & \(y \assign y- (\Delta t/2) \nabla V(x)\)\\
02 & for \(\ell = 1\) to \(L-1\)\\
03 &\qquad \(x \assign x+\Delta t M^{-1} y\)\\
04 &\qquad \(y \assign y- \Delta t \nabla V(x)\)\\
05 & end for\\
06 &\(x \assign x+\Delta t M^{-1} y\)\\
07 &\(y \assign y- (\Delta t/2) \nabla V(x)\)\\
08 & \((x^\prime, y^\prime) = (x,y)\)
\end{tabular}
\end{center}
\caption{Velocity Verlet integrator for advancing \(L\geq 1\) time-steps of length \(\Delta t>0\) starting from the initial point \(z\) and producing the final point \(z^\prime=\opInt(z)\).}
\label{fig:verlet}
\end{figure}

\begin{figure}
\begin{center}
\ttfamily\bfseries
\begin{tabular}{cll}
00 & \multicolumn{2}{l}{given \(z_0=(x_0,y_0)\in \Rd\times\Rd\) and an integer \(N\geq 1\)} \\
01 & \multicolumn{2}{l}{for \(n = 0\) to \(N-1\)}\\
02 &\qquad \(\bz_{n+1} = \oprand(z_n)\) &\% refresh momentum\\
03 &\qquad \(z_{n+1} = \opDyn(\bz_{n+1})\) &\% dynamics and accept/reject\\
04 & \multicolumn{2}{l}{end for}
\end{tabular}
\end{center}
\caption{Extra Chance Generalized HMC. It generates a Markov chain \(z_0\to z_1\to \dots \to z_N\) with \(\rho\) in \eqref{eq:distribution} as an invariant density. If \(y_0\) is drawn from the distribution \(\distrN(0,M)\),
the marginal chain of positions \(x_0\to x_1\to \dots \to x_N\) satisfies detailed balance with respect to the target distribution \eqref{eq:target}.}
\label{fig:main}
\end{figure}

\begin{figure}
\begin{center}
\ttfamily\bfseries
\begin{tabular}{cll}
00 & \multicolumn{2}{l}{given \(z=(x,y)\in \Rd\times\Rd\)}\\
01 &\multicolumn{2}{l}{\(k=0\), \(\Sigma^{(0)}(z) = 0\), \(z^{(0)} = z\)}\\
02 & draw \(u\sim U(0,1)\) &\% sample from uniform distribution\\
03 & \multicolumn{2}{l}{while \(u>\Sigma^{(k)}\) and \(k <K+1\)}\\
04 & \qquad \(z^{(k+1)} = \opInt(z^{(k)})\) &\% integration of dynamics\\
05 & \multicolumn{2}{l}{\qquad \(\Sigma^{(k+1)} = \max\Big(\Sigma^{(k)},\min\big(1, \rho(z^{(k+1)})/\rho(z^{(0)}) \big) \Big)\)}\\
06 & \multicolumn{2}{l}{\qquad \(k\assign k+1\)}\\
07 & \multicolumn{2}{l}{end while}\\
08 & \multicolumn{2}{l}{if \(u\leq \Sigma^{(k)}\) then}\\
09 & \qquad \(\opDyn(z) = z^{(k)}\) &\% \(z^{(k)}\) is accepted\\
10 & \multicolumn{2}{l}{else}\\
11 & \qquad \(\opDyn(z)=\opFlip(z^{(0)})\) &\% momentum flip\\
12 & \multicolumn{2}{l}{end if}
\end{tabular}
\end{center}
\caption{Algorithm to compute \(\opDyn(z)\). The integer \(K \geq 0\) is the number of extra chances, a user-defined parameter.}
\label{fig:D}
\end{figure}

The computation of each step of the chain in Fig.~\ref{fig:main} involves two substeps, momentum refreshment \(\oprand\) and dynamics (including an accept/reject mechanism) \(\opDyn\). Fig.~\ref{fig:D} describes the computation of the (random) transformation \(\opDyn\). The (random) mapping \(\oprand\) is defined as follows. If \(z=(x,y)\), to evaluate \(\oprand(z)\) we draw a realization \(\zeta\) from the distribution \(\distrN(0,M)\) (see \eqref{eq:distribution}) and set
\begin{equation}\label{eq:r}
\oprand (z) = (x, \cos \psi\, y+\sin\psi\, \zeta),
\end{equation}
where \(\psi\) is a user specified parameter \(0<\psi\leq \pi/2\). When \(\psi =\pi/2\) the effect of \(\oprand\) is to replace the old momentum with a fresh sample of the marginal distribution of \(y\). At the other end of the \(\psi\)-interval, choosing \(\psi \approx 0\) leads to the momentum in \(\oprand(z)\) being a small random perturbation of \(y\) as suggested in \cite{Ho91,KePe01}.

The overall algorithm in Fig.~\ref{fig:main} has the mass matrix \(M\) and the numbers \(L\), \(\Delta t\), \(\psi\) and \(K\) as parameters whose values have to be set by the user. In fact, \(M\), \(L\), \(\Delta t\), and \(\psi\) play the same role in XCGHMC as they do in generalized HMC (GHMC); it is not within the scope of the present contribution to discuss the difficult issue of how best to make the corresponding choices. Relevant references are \cite{GiCa11,BePiSa11} for the choice of \(M\), \cite{BePiRo13} for \(\Delta t\), \cite{HoGe14} for \(L\), and \cite{Ho91} for \(\psi\). Let us then study the role of \(K\), the number of \lq extra chances.\rq

As we show presently, with the choice \(K= 0\) (no extra chance) the algorithm in Fig.~\ref{fig:main} reduces to GHMC \cite{Ho91,KePe01} and therefore to standard HMC if in addition \(\psi=\pi/2\) (complete momentum refreshment). In fact, when \(K=0\), the loop \codeline{03--07} in Fig.~\ref{fig:D} is executed once per transition \(n \to n+1\) of the Markov chain. A numerical integration is performed to find the proposed state \(z^{(1)}= \opInt(\bz_{n+1}) \) and the next location \(z_{n+1}\) of the chain coincides with the proposal \(\opInt(\bz_{n+1})\) (acceptance) with probability
\[
\Sigma^{(1)} = \min\left( 1, \frac{\rho(\opInt(\bz_{n+1}))}{\rho(\bz_{n+1})}\right).
\]
In the case of rejection, the next location is \(z_{n+1} =\opFlip(\bz_{n+1})\). In both cases, acceptance and rejection, the outcome of XCGHMC coincides with that of GHMC. Regardless of the value of \(\psi\), rejection implies that \(x_{n+1} = x_n\) and this contributes to an increase of the correlation of the samples; see in this connection the results in \cite{Pe73} and \cite{Ti98}. Accordingly, rejections, while necessary for the algorithm to sample from the right distribution, are not welcome. Furthermore, consider the case where \(\psi\) has been chosen small in order that \(\oprand\) almost preserves momentum with the hope that the trajectory \(z_0, z_1, \dots z_N\) approximates the Hamiltonian dynamics (such a choice of \(\psi\) is of interest in molecular dynamics simulations). In that case, the fact that a rejection implies reversing the momentum is particularly disturbing, as it clashes with the rationale for the choice of \(\psi\).

Let us now examine the choice \(K = 1\) (one extra chance). Consider again the transition \(n\to n+1\) and assume that for given \(\bz_{n+1}\) and \(u\) (see Fig.~\ref{fig:D}), the proposal \(\opInt(\bz_{n+1})\) would have been accepted when using the standard GHMC algorithm. Under this assumption, XCGHMC goes only \emph{once} through the loop \codeline{03--07} in Fig.~\ref{fig:D} and sets \(z_{n+1} = \opInt(\bz_{n+1})\): both the outcome and the work required by XCGHM coincide with those of GHMC. On the other hand, if for given \(\bz_{n+1}\) and \(u\), the proposal \(\opInt(\bz_{n+1})\) would have been rejected in the standard GHMC algorithm, XCGHM offers \emph{one extra chance} of avoiding a momentum flip in the transition \(n\to n+1\). In fact the loop \codeline{03--07} will be executed a second time, which requires an additional integration to find \(z^{(2)}\) starting from the point \(z^{(1)}=\opInt(\bz_{n+1})\), and, if \(\Sigma^{(2)} \geq u\), the next location will be \(z_{n+1} = \opInt( \opInt(\bz_{n+1}))\) rather than \(\opFlip(\bz_{n+1})\).

The general case \(K \geq 1\) should now be clear. In Fig.~\ref{fig:D}, once a candidate \(z^{(k)}\), \(k= 1,\dots, K\), has been computed, the integration in line \codeline{04} needed to find the next candidate \(z^{(k+1)}\) will only be performed if \(z^{(k)}\) has not been accepted. A momentum flip will only occur after all candidates \(z^{(k)}\) , \(k= 1,\dots, K+1\), have failed to be accepted. Thus, with \(K \geq 1\), XCGHMC works more than GHMC in \emph{some} transitions \(n\to n+1\) in an attempt to hopefully avoid the unwelcome outcome \(z_{n+1} = \opFlip(\bz_{n+1})\). Since the bulk of the work in the algorithm lies in computing \(\opInt\), a step of the chain where the candidate \(z^{(k)}\) , \(k= 2,\dots, K+1\), is accepted costs \(k\) times as much as a step in GHMC. Numerical experiments below show that the extra cost per chain step more than pays in terms of sample quality.

The paper \cite{SoMuWe14} uses the words \lq look ahead\rq\ to describe its strategy. We feel this terminology may be misleading as it suggests that in order to accept the proposal \(\opInt(\bz_{n+1})\) it is necessary to look at the points ahead \(\opInt^2(\bz_{n+1}) = \opInt(\opInt(\bz_{n+1}))\), \(\opInt^3(\bz_{n+1}) = \opInt(\opInt(\opInt(\bz_{n+1})))\), \dots\ resulting from taking \(2L\), \(3L\), \dots\ leapfrog timesteps from the initial point \(\bz_{n+1}\). The terminology \lq delayed rejection\rq\ used in the statistical literature \cite{TiMi99,Mi01,GrMi01} is also somewhat infelicitous: the goal is to \emph{avoid} rejections, not to make them at a later time.

\section{The acceptance probabilities} \label{sec:acceptance}
Let us further study the accept/reject strategy in \(\opDyn\). Fix \(z\) and define (cf.\ line \codeline{05} in Fig.~\ref{fig:D}) a monotonic sequence
\[
0 = \Sigma^{(0)}(z) < \Sigma^{(1)}(z) \leq \dots \leq \Sigma^{(K+1)}(z) \leq 1
\]
as follows:
\begin{equation} \label{eq:xhmc-acc-rec}
\Sigma^{(0)}(z) = 0,\ \Sigma^{(k)}(z) = \max\left(\Sigma^{(k-1)}(z),\min\bigg(1, \frac{\rho(\opInt^k(z))}{\rho(z)} \bigg) \right),\ k= 1,\dots,K+1.
\end{equation}
We emphasize that in this section we are concerned with the \emph{analysis} of \(\opDyn\): if \(k >1\), the actual algorithm will \emph{not} compute \(\Sigma^{(k)}(z)\) (which requires \(k\) integration legs with \(L\) time-steps each to find \(\opInt^k(z)\)) unless the point \(\opInt^{k-1}(z)\) has been rejected (as its predecessors must have).
Note that, for \(k = 1, \dots, K+1\),
\begin{equation} \label{eq:xhmc-acc}
\Sigma^{(k)}(z)
= \max_{1\leq j \leq k} \min\left(1, \frac{\rho(\opInt^j(z))}{\rho(z)}\right)
= \min \left(1, \max_{1\leq j \leq k}\frac{\rho(\opInt^j(z))}{\rho(z)}\right).
\end{equation}

From Fig.~\ref{fig:D} it is easily concluded that, when computing the random point \(\opDyn(z)\), the event \lq \(\opDyn(z)\) is one of \(\opInt(z)\), \(\opInt^2(z)\), \dots, \(\opInt^k(z)\)\rq\ happens if and only if \(u \leq \Sigma^{(k)}(z)\); therefore that event has probability \(\Sigma^{(k)}(z)\). As a consequence, the difference
\begin{equation} \label{eq:xhmc-trns}
p^{(k)}(z) = \Sigma^{(k)}(z)-\Sigma^{(k-1)}(z), \qquad k = 1,\dots,K+1
\end{equation}
represents the probability that \(\opDyn(z) = \opInt^k(z)\) and
\begin{equation}\label{eq:pk}
p^{(K+2)}(z) = 1- \Sigma^{(K+1)}(z)
\end{equation}
is the probability that none among \(\opInt(z)\), \(\opInt^2(z)\), \dots, \(\opInt^{K+1}(z)\) is accepted and therefore \(\opDyn(z) = \opFlip(z)\). Pictorially, the values \(\Sigma^{(k)}(z)\), \(k= 1,\dots, K+1\), partition the interval \([0,1]\) into \(K+2\) subintervals whose lengths \(p^{k}(z)\), \(k= 1,\dots, K+2\), provide the probabilities of the outcomes \(\opDyn(z) = \opInt(z)\), \dots, \(\opDyn(z)= \opInt^{K+1}(z)\), \(\opDyn(z) = \opFlip(z)\). Some of the subintervals may degenerate into a single point and then the corresponding outcomes have zero probability.

We note the following conclusions:
\begin{enumerate}
\item Assume that \(1 <k\leq K+1\) and \(\rho(\opInt^k(z))\) does not exceed the maximum of the preceding \(\rho(\opInt^j(z))\), \(j = 1,\dots, k-1\). Then \(\Sigma^{(k)}(z)=\Sigma^{(k-1)}(z)\) and the event \(\opDyn(z) = \opInt^k(z)\) has probability \(p^{(k)}(z) =0\): \emph{\(\opDyn\) only moves the chain from \(z\) to locations \(\opInt^k(z)\) where the density \(\rho\) is larger than at all the preceding \lq missed chances\rq\ \(\opInt^j(z)\), \(j = 1,\dots, k-1\).}
\item Assume that the value of \(\rho\) at one location \(\opInt^j(z)\), \(1 \leq j \leq k\), exceeds the value \(\rho(z)\). Then \(\Sigma^{(k)}(z) = 1\) and accordingly one among \(\opInt(z)\), \(\opInt^2(z)\), \dots, \(\opInt^k(z)\) will be accepted. In particular, { \em if \(\rho(\opInt^k(z)) > \rho(z)\) for at least a value of \(k\), \(1 \leq k \leq K+1\), the momentum flip \(\opDyn(z) = \opFlip(z)\) will not occur.}
\end{enumerate}

It is well known that for \emph{symplectic integrators} \cite{SaCa94,HaLuWa06,LeRe04}, including the leapfrog scheme, the value of \(H\) along a long numerical trajectory of a Hamiltonian system typically oscillates around the value at the initial point of the integration. Here (see \eqref{eq:distribution}) this means that, if \(K\) is large, some of the values \(\rho(\opInt^k(z))\) may be expected to be below \(\rho(z)\). In view of item 2 above, it then may be hoped that then the flip \(\opDyn(z) = \opFlip(z)\) will not occur (see the numerical experiments below).

When comparing the present approach with the literature on delayed rejection \cite{TiMi99,Mi01,GrMi01}, it is useful to emphasize that, if \(z\) is given, the probabilities \(\Sigma^{(k)}(z)\) and \(p^{(k)}(z)\) defined above are not conditional. We have not been concerned with questions such as what is the probability of accepting \(\opInt^2(z)\) \emph{conditional} on \(\opInt(z)\) having already been rejected. This is to be compared with the approach in \cite{TiMi99,Mi01,GrMi01}, which focuses on the probability of accepting a new attempt conditional on previous attempts at the same step having been rejected. The formulas for such conditional probabilities turn out not to be very neat.

\section{Numerical results} \label{sec:numerics}
We have implemented the algorithm in Fig.~\ref{fig:main} for sampling from the canonical distribution of a molecule of \(C_9H_{20}\). This example, taken from \cite{RyBe78}, has been used in \cite{CaLeSt07} to compare different Markov chain samplers. It has \(3\times 9 = 27\) degrees of freedom as the hydrogen atoms are lumped to the corresponding carbon atoms. The degrees of freedom may of course be associated with the six rigid body motions and with vibrations in the eight carbon-carbon bond lengths, seven bond angles and six dihedral angles. The potential involves two-, three- and four-atom contributions related to bond lengths, bond angles and dihedral angles respectively and furthermore Lenard-Jones terms for all pairs of atoms separated by three or more covalent bounds. Units and parameter values here are as in \cite{CaLeSt07}, with the inverse temperature \(\beta = 1\).

The potential energy associated with each dihedral angle has three possible stable equilibria. We are interested in estimating the probability that the absolute value of the dihedral angle \(\phi_1\) between the first four atoms is below 1.75; this corresponds with \(\phi_1\) lying in the basin of attraction of the most stable value \(\phi_1 = 0\) (for which the first four carbon atoms are coplanar). In other words we assume that the samples \(x_n\) are to be used in \eqref{eq:average} when \(A\) is the indicator of the set \(\{|\phi_1| \leq 1.75\}\), i.e.\ \(A = 1\) when \(|\phi_1| \leq 1\) and \(A=0\) if \(|\phi_1| > 1\). The quality of the samples generated by the different algorithms will be measured by the \emph{effective sample size} (relative to to this choice of \(A\)); the ESS of a set of \(N+1\) (correlated) Markov chain samples \(x_n\) represents the number of \emph{independent} samples that contain the same amount of information; in other words the ESS is the result of dividing \(N+1\) by the autocorrelation time. We estimate ESSs by means of the initial monotone estimator in \cite{Ge92}. Several alternative choices of the observable \(A\) were also considered; while ESSs change substantially with \(A\), the corresponding numerical results lead to the same qualitative conclusions as those reported here.

All experiments have a unit mass matrix as in \cite{CaLeSt07}. We studied the cases \(K=0\) (GHMC) and \(K=3\) (three extra chances); the results below show very clearly that it is unnecessary to consider higher values of \(K\). Dozens of combinations of the remaining parameters \(\Delta t\), \(L\) and \(\psi\) were tried and we only report here on a representative selection. For each choice of parameters in the algorithm, we generated ten realizations of the chain. The values of the ESS and acceptance probability given here are averages over the ten realizations. Each realization included an initial burn-in phase (500 transitions) and a production phase; the latter comprised \(10^6\) evaluations of the force \(-\nabla V\) so that all realizations use the same amount of computational work, regardless of the value of \(L\). In other words, increasing \(L\) results in a lower number \(N\) of computed samples. Note also that when \(K=3\) a higher use of extra chances in a realization leads to fewer samples being generated. To complete the description of our experimental setting we mention that, in order to avoid resonances, the value of \(\Delta t\) was slightly randomized at the beginning of each integration leg \(\opInt\) by subjecting it to a perturbation ranging in a \(\pm 5\%\) interval, see e.g. \cite{Ne11}.

\begin{table}
\begin{center}
\begin{tabular}{c||c||cccc|c}
\(\Delta t\) &
\(K=0\) &
\multicolumn{5}{c}{\(K=3\)}\\
\hline
0.012 & 93\% & 93\% &  6\% & 1\% & 0\% & 100\%\\
0.016 & 86\% & 87\% & 11\% & 2\% & 0\% & 100\%\\
0.020 & 77\% & 80\% & 16\% & 3\% & 1\% & 100\%\\
0.024 & 65\% & 71\% & 22\% & 6\% & 1\% & 100\%\\
\multicolumn{2}{c}{}&\(\scriptstyle a_0\) & \(\scriptstyle a_1\) & \(\scriptstyle a_2\) & \multicolumn{1}{c}{\(\scriptstyle a_3\)} & \({\scriptscriptstyle \sum}\scriptstyle  a_k\)
\end{tabular}
\end{center}
\caption{Acceptance rates (rounded to the nearest percentage point) when \(\sin \psi = 1\) and \(L\Delta t = 0.48\) for different values of \(\Delta t\) and \(K=0\) (standard HMC) and \(K=3\) (three extra chances). The column labelled \(a_k\) gives the number corresponding to acceptance taking place after \(k\) extra chances.}
\label{tab:acceptances}
\end{table}
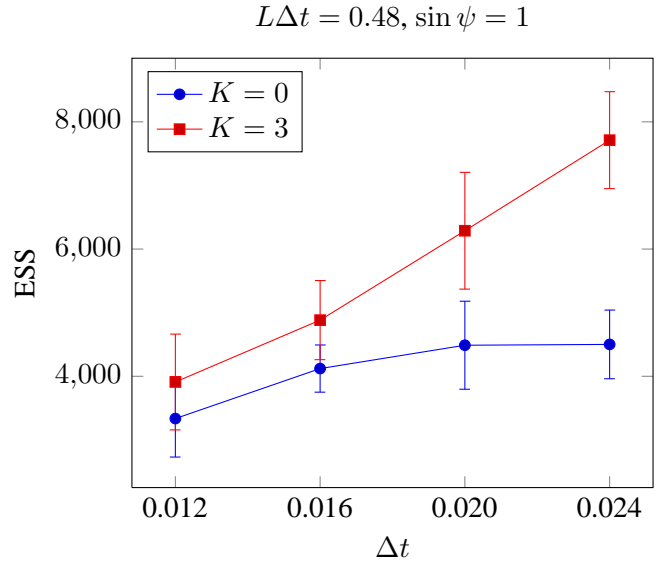
\begin{figure}
\begin{center}
\begin{tikzpicture}
  \begin{axis}[
    xlabel = {\(\Delta t\)},
    ylabel = {ESS},
    legend entries = {\(K=0\),\(K=3\)},
    legend pos = {north west},
    title = {\(L\Delta t=0.48\), \(\sin\psi=1\)},
    ymin =  2250, ymax =  9000,
    xtick = {0.012,0.016,0.020,0.024},
    xticklabels = {0.012,0.016,0.020,0.024},
    scaled x ticks = false,
    ]
    \addplot+[error bars/.cd,y dir=both,y explicit] coordinates {
      (0.012,3335) +- (0,605)
      (0.016,4121) +- (0,371)
      (0.020,4488) +- (0,692)
      (0.024,4501) +- (0,540)
    };
    \addplot+[error bars/.cd,y dir=both,y explicit] coordinates {
      (0.012,3910) +- (0,753)
      (0.016,4883) +- (0,622)
      (0.020,6288) +- (0,918)
      (0.024,7712) +- (0,762)
    };
  \end{axis}
\end{tikzpicture}
\end{center}
\caption{Effective sample size as a function of \(\Delta t\) when the integration interval is \(L\Delta t = 0.48\) and \(\sin\psi = 1\) (complete momentum refreshment).}
\label{fig:delta}
\end{figure}
\begin{figure}
\begin{center}
\begin{tikzpicture}
  \begin{axis}[
    xlabel = {\(L\Delta t\)},
    ylabel = {ESS},
    legend entries = {\(K=0\),\(K=3\)},
    legend pos = {north west},
    title = {\(\Delta t=0.024\), \(\sin\psi=1\)},
    ymin = 2250, ymax = 9000,
    xtick = {0.12,0.24,0.36,0.48,0.60,0.72},
    xticklabels = {0.12,0.24,0.36,0.48,0.60,0.72},
    ]
    \addplot+[error bars/.cd,y dir=both,y explicit] coordinates {
      (0.12,2675) +- (0,697)
      (0.24,3942) +- (0,875)
      (0.36,3558) +- (0,587)
      (0.48,4501) +- (0,540)
      (0.60,4844) +- (0,948)
      (0.72,3970) +- (0,511)
    };
    \addplot+[error bars/.cd,y dir=both,y explicit] coordinates {
      (0.12,3963) +- (0,650)
      (0.24,6129) +- (0,793)
      (0.36,6596) +- (0,479)
      (0.48,7712) +- (0,762)
      (0.60,7438) +- (0,617)
      (0.72,7189) +- (0,1259)
    };
  \end{axis}
\end{tikzpicture}
\end{center}
\caption{Effective sample size as a function of the time-span \(L\Delta t\) of the numerical integration when the time-step is \(\Delta t = 0.024\) and \(\sin\psi = 1\) (complete momentum refreshment).}
\label{fig:L}
\end{figure}
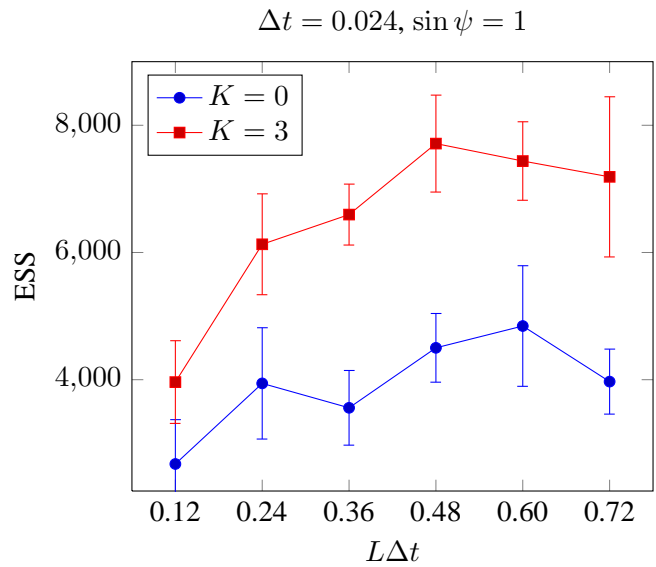
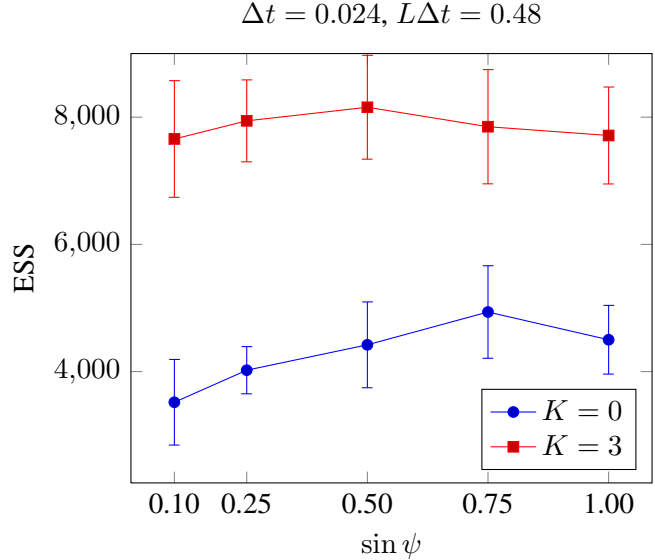
\begin{figure}
\begin{center}
\begin{tikzpicture}
  \begin{axis}[
    xlabel = {\(\sin\psi\)},
    ylabel = {ESS},
    legend entries = {\(K=0\),\(K=3\)},
    legend pos = {south east},
    title = {\(\Delta t=0.024\), \(L\Delta t=0.48\)},
    ymin = 2250, ymax = 9000,
    xtick = {0.10,0.25,0.50,0.75,1.00},
    xticklabels = {0.10,0.25,0.50,0.75,1.00},
    scaled x ticks = false,
    ]
    \addplot+[error bars/.cd,y dir=both,y explicit] coordinates {
      (0.10,3518) +- (0,673)
      (0.25,4022) +- (0,371)
      (0.50,4421) +- (0,675)
      (0.75,4937) +- (0,727)
      (1.00,4501) +- (0,540)
    };
    \addplot+[error bars/.cd,y dir=both,y explicit] coordinates {
      (0.10,7657) +- (0,917)
      (0.25,7942) +- (0,644)
      (0.50,8156) +- (0,816)
      (0.75,7850) +- (0,898)
      (1.00,7712) +- (0,762)
    };
  \end{axis}
\end{tikzpicture}
\end{center}
\caption{Effective sample size as a function of \(\sin\psi\) when the time-step is \(\Delta t = 0.024\) and \(L\Delta t = 0.48\).}
\label{fig:psi}
\end{figure}

Table~\ref{tab:acceptances} gives values of the acceptance probability when \(L\Delta t = 0.48\) and \(\sin \psi = 1\) with different values of \(\Delta t\). The maximum value \(\Delta t = 0.024\) considered is fairly close to the upper limit allowed by the size of the stability interval of the Verlet integrator, which experiments suggest is \(\approx 0.030\). As expected, for \(K=0\) (standard HMC) the fraction of accepted steps decreases as \(\Delta t\) increases (lower accuracy in the integration, higher energy errors). It is apparent that the use of extra chances almost completely eliminates rejections (in the last row, with two decimal places, \(\sum a_k=99.80\%\)). For small \(\Delta t\) the value of \(a_0\) coincides with the acceptance rate for \(K=0\) and in fact, it is clear from our discussion of the algorithm that, when the chains are at \emph{stationarity}, the \emph{expected value} of \(a_0\) exactly matches the expected acceptance rate of HMC. The agreement between the values of \(a_0\) and the acceptance rate for \(K=0\) deteriorates as \(\Delta t\) increases; this happens because for large \(\Delta t\) there are wider differences between the empirical values obtained from the simulations and the corresponding theoretical expected values at stationarity. For large \(\Delta t\), XGCHMC is less prone to be stuck at the present location than HMC. 

Fig.~\ref{fig:delta} corresponds to the parameter values in Table~\ref{tab:acceptances} and gives the ESS of the samples generated. Comparing the values \(\Delta t = 0.020\) and \(\Delta t = 0.024\) for HMC (i.e.\ \(K=0\)) we see that a lower acceptance rate does not automatically imply a less effective simulation: a larger \(\Delta t\) requires less work per integration leg and provides more samples with a given amount of work. Note that an acceptace rate close to \(65\%\) in HMC is sometimes regarded as optimal (see \cite{BePiRo13} and its references). The figure clearly shows that the additional work per transition required by the extra chances pays: in \emph{all} simulations the ESS improves substantially when moving from \(K=0\) to \(K=3\). In this figure, the best ESS for the extra chance algorithm (7712) is more than \(70\%\) higher the best ESS for the its standard counterpart (4501).

In Fig.~\ref{fig:L} we study the variation of ESS as a function of the time-span \(L\Delta t\) of the numerical integration when \(\Delta t = 0.024\) and \(\sin \psi = 1\). For \(K=0\), ESS presents a clear dip at \(L\Delta t = 0.36\), probably due to a resonance between the integration time-span and the periods involved in the dynamics of \(\phi_1\). For all values of \(L\Delta t\), the ESS with extra chances improves clearly on that of the standard GHMC. The acceptance rates for the simulations in this figure do not differ much from those given in the last row of Table~\ref{tab:acceptances}: for instance \(L\Delta t = 0.12\) has \(67\%\) for \(K=0\) and \(74\%+19\%+5\%+1\%\) for \(K=3\).

Finally Fig.~\ref{fig:psi} shows the dependence of ESS on Horowitz's angle \(\psi\). As it may have been expected, extra chances are more beneficial for small values of \(\psi\). Again the acceptance rates here are roughly the same as those in the last row of Table~\ref{tab:acceptances}. We conclude that for fixed \(\Delta t\) the acceptance rates are almost independent of the time-span \(L\Delta t\) and \(\psi\).

\section{Conclusion} \label{sec:conclusion}
We have studied a technique to avoid rejections/momentum flips in the HMC and GHMC methods. Proof of concept experiments show clearly that the technique is very promising in terms of the number of uncorrelated samples that may be obtained with a given amount of computational work.

The material here may be extended in different directions. It may be combined with more sophisticated integrators \cite{BlCaSa14,PrLiEa12} (particularly so if the dimensionality of the target is very high and HMC requires small energy errors per degree of freedom) or with the use of shadow Hamiltonians \cite{IzHa04,AkRe08}. It may also be easily extended to cover \lq compressible\rq\ variants of HMC such as those presented in \cite{FaSaSk14}.

\section*{Acknowledgments}
We are grateful to M. P. Calvo, M. Girolami, T. Radivojevic, R. D. Skeel and J. Sohl-Dickstein for their inputs and to the ICMAT for the usage of their cluster. This research is supported by Ministerio de Ciencia e Innovación (Spain), under the project MTM2010-18246-C03-01, and by Junta de Castilla y León (Spain) together with the European Social Fund, through a postdoctoral position.

\appendix
\renewcommand{\thesection}{\Alph{section}}
\section{Stationarity of the Markov chain} \label{sec:stationarity}

In this appendix we prove that the distribution \eqref{eq:distribution} is invariant for the chain in Fig.~\ref{fig:main}. Since \(\rho\) clearly remains invariant under \(\oprand\), it is sufficient to show that it also remains invariant under \(\opDyn\). It is well known (see e.g.\ \cite{LeRoSt10,FaSaSk14,Sa14}) that, in turn, such an invariance is implied by the requirement of \emph{modified detailed balance}: for each \(z\) and \(z^\prime\),
\begin{equation}
\label{eq:mdb}
\rho(z)\, \rho(z^\prime\mid z) = \rho(\opFlip(z^\prime))\, \rho(\opFlip(z)\mid \opFlip(z^\prime)).
\end{equation}
Here \(\rho(z^\prime \mid z)\) denotes the conditional distribution of \(z^\prime=\opDyn(z)\) conditional on \(z\) and \(\rho(\opFlip(z)\mid \opFlip(z^\prime))\) the conditional distribution of \(\opFlip(z) = \opDyn(\opFlip(z^\prime))\) conditional on \(\opFlip(z^\prime)\) (i.e. on \(z^\prime\)).

In order to establish \eqref{eq:mdb}, we shall employ the following identity, which is valid for each \(z\) in phase space and \(k= 1,\dots,K+1\) and will be proved at the end of this appendix,
\begin{equation}\label{eq:main}
\rho(z)\, p^{(k)}(z) = \rho(\opFlip(\opInt^k(z)))\,p^{(k)}(\opFlip(\opInt^k(z))).
\end{equation}

By definition of \(\opDyn\), the left and right hand-sides of \eqref{eq:mdb} have the values
\[
\rho(z) \left(\sum_{k=1}^{K+1} p^{(k)}(z)\, \delta(z^\prime-\opInt^k(z)) + p^{(K+2)}(z) \,\delta(z^\prime-\opFlip(z))
\right)\]
and
\[
\rho(\opFlip(z^\prime)) \left(\sum_{k=1}^{K+1} p^{(k)}(\opFlip(z^\prime))\, \delta(\opFlip(z)-\opInt^k(\opFlip(z^\prime))) + p^{(K+2)}(\opFlip(z^\prime))\, \delta( \opFlip(z)-\opFlip(\opFlip(z^\prime)) )\right)
\]
respectively. Since the last terms in these expressions obviously coincide, the proof will be over if we show that, for \(k= 1,\dots, K+1\),
\[
\rho(z)\, p^{(k)}(z)\, \delta(z^\prime-\opInt^k(z)) = \rho(\opFlip(z^\prime))\, p^{(k)}(\opFlip(z^\prime))\, \delta(\opFlip(z)-\opInt^k(\opFlip(z^\prime))).
\]
This is a consequence of \eqref{eq:main}, as the reversibility of \(\opInt\) implies
\[
\delta(\opFlip(z)-\opInt^k(\opFlip(z^\prime))) = \delta (z- \opFlip(\opInt^k(\opFlip(z^\prime)))) = \delta(z^\prime-\opInt^k(z))
\]
because both \(\opFlip\) and \(\opInt^k\) are volume preserving. The proof of \eqref{eq:mdb} is then ready.

\begin{figure}
\[\begin{array}{ccccccccc}
z \phantom{\opFlip} & \to & \dots & \to & \opInt^j(z)& \to & \dots & \to & \opInt^k(z)\vspace{1ex}\\
\updownarrow \opFlip &&&&\updownarrow \opFlip&&&&\updownarrow \opFlip\vspace{1ex}\\
\opInt^k(\opFlip(\opInt^k(z))) & \gets & \dots & \gets & \opInt^{k-j}(\opFlip(\opInt^k(z)))& \gets & \dots & \gets & \opFlip(\opInt^k(z))
\end{array}\]
\caption{Top row: \(k\) successive applications of \(\opInt\) starting from the point \(z\). Bottom row: \(k\) successive applications of \(\opInt\) starting from \(\opFlip(\opInt^k(z))\). The reversibility of \(\opInt\) implies that each point in the bottom row may be obtained from the corresponding point in the top row by flipping the momentum. Therefore each point in the bottom row possesses the same density \(\rho\) as the corresponding point above.}
\label{fig:reversibility}
\end{figure}

Let us finally prove the identity \eqref{eq:main}. We restrict the attention to the case \(k>1\); the proof for the case \(k=1\) is similar but simpler. After using the expression for \(p^{(k)}\) in \eqref{eq:pk}, we have to show that
\begin{multline*}
\min \left(
\rho(z),
\max_{1\leq j \leq k}\rho(\opInt^j(z))\right)-
\min \left(
\rho(z),
\max_{1 \leq j \leq k-1}\rho(\opInt^j(z))\right)
\\=
\min \left(
\rho(\opFlip(\opInt^k(z))),
\max_{1\leq j \leq k}\rho(\opInt^j( \opFlip(\opInt^k(z))))\right)-
\min \left(
\rho(\opFlip(\opInt^k(z))),
\max_{1 \leq j \leq k-1}\rho(\opInt^j( \opFlip(\opInt^k(z))))\right).
\end{multline*}
The reversibility of \(\opInt\) implies (see Fig.~\ref{fig:reversibility}) that the values of \(\rho\) that feature on the right hand-side are actually the same as those appearing in the left hand side. In fact, if \(\hj\), \(1\leq \hj\leq k-1\), is such that
\[
\rho(\opInt^\hj(z)) = \max_{1 \leq j \leq k-1}\rho(\opInt^j(z)),
\]
the equality to be established reads
\begin{multline*}
\min \left(
\rho(z),
\max\left( \rho(\opInt^\hj(z)), \rho(\opInt^k(z))\right)\right)-
\min \left(
\rho(z),
\rho(\opInt^\hj(z))\right)
\\=
\min \left(
\rho(\opInt^k(z)),
\max\left( \rho(z),\rho(\opInt^\hj(z))\right)\right)-
\min \left(
\rho(\opInt^k(z)),\rho(\opInt^\hj(z))\right).
\end{multline*}
This is checked by successively considering the six possible orderings \(\rho(z) \leq \rho(\opInt^\hj(z)) \leq \rho(\opInt^k(z))\), \(\rho(z) \leq \rho(\opInt^k(z)) \leq \rho(\opInt^\hj(z))\), etc.

\section{Detailed balance} \label{sec:balance}
\begin{figure}
\begin{center}
\ttfamily\bfseries
\begin{tabular}{cll}
00 & \multicolumn{2}{l}{given \(Z_0=(X_0,Y_0)\in \Rd\times\Rd\) and an integer \(N\geq 1\)}\\
01 & \multicolumn{2}{l}{for \(n = 0\) to \(N-1\)}\\
02 & \qquad \(\bZ_{n+1} = \opRand(Z_n)\) &\% refresh momentum\\
03 & \qquad \(\hZ_{n+1} = \opDyn(\bZ_{n+1})\) &\% dynamics and accept/reject\\
04 & \qquad \(Z_{n+1} = \opRand(\hZ_{n+1})\) &\% refresh momentum\\
05 & \multicolumn{2}{l}{end for}
\end{tabular}
\end{center}
\caption{A \lq tought experiment\rq\ chain that satisfies modified detailed balance with respect to \(\rho\) in \eqref{eq:distribution}.}
\label{fig:tought}
\end{figure}

While, as proved above, a single step of \(\opDyn\) satisfies modified detailed balance, the same is not true for the whole chain in Fig.~\ref{fig:main}. The paper \cite{FaSaSk14} presents a relevant counterexample and also shows that modified detailed balance is achieved by imposing a palindromic structure to the transitions \(n\to n+1\) as in Fig.~\ref{fig:tought}. The lines \codeline{02} and \codeline{04} refresh the momentum as in \eqref{eq:r}; the use of a capital \(\opRand\) indicates that Fig.~\ref{fig:tought} may employ an angle \(\Psi\), \(0<\Psi\leq \pi/2\) different from that \(\psi\) in Fig.~\ref{fig:main}:
\[
\opRand (Z) = (X, \cos \Psi\, Y+\sin\Psi\, \zeta).
\]
In fact in what follows we assume that \(\Psi\) is determined as a function of \(\psi\) through the relation \(\cos ^2 \Psi = \cos \psi\).

We show next that the XCGHMC chain (Fig.~\ref{fig:main}) implemented in practice and the \lq thought experiment chain\rq\ in Fig.~\ref{fig:tought} give rise to same marginal chain for the positions \(x_n\) or \(X_n\) and that, therefore, the chain of positions in XCGHMC is reversible with respect to the target \eqref{eq:target}.

Let us consider an arbitrary realization of the auxiliary chain in Fig.~\ref{fig:tought}. This will correspond to a starting location \((X_0,Y_0)\) (with \(Y_0\) drawn from the distribution \(\distrN(0,M)\)) and to realizations \(\bzeta_{n+1}\), \(u_{n+1}\), \(\hzeta_{n+1}\) of the random variables used in steps \codeline{02}, \codeline{03} and \codeline{04} respectively. We define a realization of the XCGHMC chain by setting \(x_0= X_0\),
\begin{eqnarray*}
y_0 &=& \cos(\psi-\Psi)\, Y_0-\sin(\psi-\Psi)\bzeta_1,\\
\zeta_1 &= & \sin (\psi-\Psi)\, Y_0+ \cos(\psi-\Psi)\bzeta_1,\\
\zeta_{n+1} & = & (1/\sin \psi) (\cos \Psi\, \sin\Psi \hzeta_n+ \sin \Psi \bzeta_{n+1}),\quad n = 1,\dots, N-1,
\end{eqnarray*}
and using \(\zeta_{n+1}\) and \(u_{n+1}\), \(n = 0,\dots,N-1\), when computing \(\oprand\) and \(\opDyn\) in steps \codeline{02} and \codeline{03} of Fig.~\ref{fig:main} respectively. Note that this is legitimate because, as it is easily checked, the random variables \(\zeta_n\) possess the right \(\distrN(0,M)\) distribution. The definitions of \(y_0\) and \(\zeta_1\) result in \(\bz_1 = \bZ_1\), which in turn implies that \(z_1 = \hZ_1\). Then the choice of \(\zeta_2\) ensures that the momentum refreshment to get \(\bz_2 = \oprand(z_1)\) produces the same output as the two momentum refreshments to compute \(\bZ_2 = \opRand(Z_1) = \opRand(\opRand(\hZ_1))\). The iteration of this argument shows that, for \(n= 1,\dots, N\), \(z_n =\hZ_n\) so that \(x_n= X_n\). We conclude that we may regard the samples \(x_n\) as originating from the palindromic chain in Fig.~\ref{fig:tought}.

\section{Equivalence of XCGHMC and the algorithm by Sohl-Dickstein \emph{et al.}}  \label{sec:equivalence}

In \cite{SoMuWe14}, the authors introduced an algorithm called Look Ahead Hybrid Monte Carlo (LAHMC) which is in fact equivalent to XCGHMC. There, the transition probabilities to \(z^{(k)}\) are defined as follows 
\begin{equation} \label{eq:lahmc-trns}
\pi^{(k)}(z) = \min\Bigg(
1-\sum_{1\leq j\leq k-1}\pi^{(j)}(z),
\frac{\rho(\opFlip\opInt^k(z))}{\rho(z)}\bigg(1-\sum_{1\leq j\leq k-1}\pi^{(j)}(\opFlip\opInt^k(z))\bigg)
\Bigg), \ \ k=1,\dots,K+1.
\end{equation}

In order to show that XCGHMC and LAHMC are equivalent, it suffices to show that the accumulated probabilities
\begin{equation} \label{eq:lahmc-acc}
S^{(0)}(z) = 0,\ \
S^{(k)}(z) = \sum_{1\leq j\leq k}\pi^{(j)}(z),\ \
k=1,\dots,K+1.
\end{equation}
agree with \(\Sigma^{(k)}\) in \eqref{eq:xhmc-acc} for all \(k\) at any point \(z\). Simple manipulations prove that this is true for \(k=1\); the general case is proved by induction. Assume the assertion true up to some \(1\leq k-1\leq K\) and let us show that it holds for \(k\). To lighten the writing, we set
\[ M^{(k)}(z) = \max_{1\leq j \leq k}\frac{\rho(\opFlip\opInt^j(z))}{\rho(z)} = \max_{1\leq j \leq k}\frac{\rho(\opInt^j(z))}{\rho(z)} . \]
By definition,
\begin{eqnarray*}
S^{(k)}(z)
&=& S^{(k-1)}(z) + \min\bigg(1-S^{(k-1)}(z),\frac{\rho(\opFlip\opInt^k(z))}{\rho(z)}\Big(1-S^{(k-1)}(\opFlip\opInt^k(z))\Big)\bigg)\\
&=& \min\bigg(1,\frac{\rho(\opFlip\opInt^k(z))}{\rho(z)}\Big(1-S^{(k-1)}(\opFlip\opInt^k(z))\Big) + S^{(k-1)}(z)\bigg)\\
&=& \min\bigg(1,\frac{\rho(\opFlip\opInt^k(z))}{\rho(z)}\Big(1-\min \big(1, M^{(k-1)}(\opFlip\opInt^k(z))\big)\Big) + S^{(k-1)}(z)\bigg)\\
&=& \min\bigg(1,\max \Big(0,\frac{\rho(\opFlip\opInt^k(z))}{\rho(z)}-M^{(k-1)}(z)\Big) + S^{(k-1)}(z)\bigg),
\end{eqnarray*}
where we have used the hypothesis of induction on \(S^{(k-1)}(\opFlip\opInt^k(z))\) at the third equality and the fact that
\[ \rho(z)\,M^{(k-1)}(z)=\rho(\opFlip\opInt^k(z))\,M^{(k-1)}(\opFlip\opInt^k(z)) \]
at the last equality (recall that \(\rho\) is invariant under \(\opFlip\) and that \(\opInt^j\opFlip\opInt^k=\opInt^{k-j}\)).

Now, two possibilities arise, either \(S^{(k-1)}(z)<1\) or \(S^{(k-1)}(z)=1\). If \(S^{(k-1)}(z)<1\), in which case by the hypothesis of induction \(S^{(k-1)}(z)=M^{(k-1)}(z)\), then
\[ S^{(k)}(z) = \min\bigg(1,\max\Big(M^{(k-1)}(z),\frac{\rho(\opFlip\opInt^k(z))}{\rho(z)}\Big)\bigg). \]
If on the contrary \(S^{(k-1)}(z)=1\), the hypothesis of induction implies that \(M^{(k-1)}(z)\geq1\), hence
\[ S^{(k)}(z) = 1 = \min\bigg(1,\max\Big(M^{(k-1)}(z),\frac{\rho(\opFlip\opInt^k(z))}{\rho(z)}\Big)\bigg). \]
Therefore, in both cases
\[ S^{(k)}(z) = \min\bigg(1,\max\Big(M^{(k-1)}(z),\frac{\rho(\opFlip\opInt^k(z))}{\rho(z)}\Big)\bigg) = \min\Big(1,M^{(k)}(z)\Big), \]
and the proof is complete.


In our opinion working with accumulated probabilities as in XCGHMC rather than with the probabilities \(\pi^{(k)}\) not only leads to better theoretical insights as in Section \ref{sec:acceptance}, but results in simpler code.

\bibliographystyle{elsarticle-harv}
\bibliography{hmc}
\end{document}